\documentclass[12pt,draftcls,onecolumn]{IEEEtran}

\usepackage{url}
\usepackage{amssymb}
\usepackage{makeidx}         
\usepackage{amsmath}
\usepackage{multicol}        
\usepackage[bottom]{footmisc}
\usepackage{graphicx}
\usepackage{color}
\usepackage{amsfonts}
\usepackage[latin1]{inputenc}
\usepackage{epstopdf}

\newcommand{\eq}{\begin{equation}}
\newcommand{\eeq}{\end{equation}}
\newcommand{\eqn}{\begin{eqnarray}}
\newcommand{\eeqn}{\end{eqnarray}}

\newcommand{\bsea}{\begin{subeqnarray}}
\newcommand{\esea}{\end{subeqnarray}}
\newcommand{\nn}{\nonumber}

\newcommand{\Sp}[2]{\left< #1,#2 \right> }

\newcommand{\Min}[1]{\,\underset{#1}{\mathrm{min}}\,}

\newcommand{\tr}{\mathop{\rm tr}}  
\newcommand{\e}[1]{\mathrm{e}^{#1}}
\newcommand{\Set}[1]{\left\{ #1\right\}}

\newcommand{\Cc}{ \mathcal{C}}
\newcommand{\Dc}{ \mathcal{D}}

\newcommand{\Fc}{ \mathcal{F}}

\newcommand{\Ic}{ \mathcal{I}}
\newcommand{\Jc}{ \mathcal{J}}
\newcommand{\Kc}{ \mathcal{K}}
\newcommand{\Lc}{ \mathcal{L}}

\newcommand{\Qc}{ \mathcal{Q}}

\newcommand{\Sc}{ \mathcal{S}}

\newcommand{\Vc}{ \mathcal{V}}

\newcommand{\Cs}{ \mathbb{C}}
\newcommand{\Ds}{ \mathbb{D}}
\newcommand{\Es}{ \mathbb{E}}

\newcommand{\Hs}{ \mathbb{H}}

\newcommand{\Ns}{ \mathbb{N}}

\newcommand{\Rs}{ \mathbb{R}}

\newcommand{\Ts}{ \mathbb{T}}

\newcommand{\Zs}{ \mathbb{Z}}

\def\qed{\hfill \vrule height 7pt width 7pt depth 0pt \smallskip}

\newcounter{pippo}

\newtheorem{remark}{Remark}[section]
\newtheorem{teor}{Theorem}[section]
\newtheorem{corr}{Corollary}[section]
\newtheorem{propo}{Proposition}[section]
\newtheorem{lemm}{Lemma}[section]
\newtheorem{exam}{Example}
\newtheorem{probl}[pippo]{Problem}
\newtheorem{defn}{Definition}[section]

\newcommand{\teo}{\begin{teor}}
\newcommand{\eteo}{\end{teor}}
\newcommand{\cor}{\begin{corr}}
\newcommand{\ecor}{\end{corr}}
\newcommand{\prop}{\begin{propo}}
\newcommand{\eprop}{\end{propo}}
\newcommand{\lem}{\begin{lemm}}
\newcommand{\elem}{\end{lemm}}
\newcommand{\ex}{\begin{exam}}
\newcommand{\eex}{\end{exam}}
\newcommand{\pb}{\begin{probl}}
\newcommand{\epb}{\end{probl}}
\newcommand{\df}{\begin{defn}}
\newcommand{\edf}{\end{defn}}
\newcommand{\aprop}{\begin{apropo}}
\newcommand{\eaprop}{\end{apropo}}
\newcommand{\alem}{\begin{alemm}}
\newcommand{\ealem}{\end{alemm}}
\newcommand{\rem}{\begin{remark}}
\newcommand{\erem}{\end{remark}}

\newcommand{\Lp} {\mathcal{L}_+}

\newcommand{\Stm} {\Sc_+^m(\Ts)}
\begin{document}

\title{Multivariate Spectral Estimation based on the concept of Optimal Prediction}

\author{Mattia Zorzi}

\author{Mattia~Zorzi\thanks{This note presents research
results of the Belgian Network DYSCO (Dynamical Systems, Control, and
Optimization), funded by the Interuniversity Attraction Poles Programme, initiated
by the Belgian State, Science Policy Office. The scientific responsibility
rests with its authors. This research is also supported by FNRS (Belgian Fund for Scientific Research).}
\thanks{M. Zorzi is  with the
Dipartimento di Ingegneria dell'Informazione, Universit\`a di
Padova, via Gradenigo 6/B, 35131 Padova, Italy, {\tt\small
zorzimat@dei.unipd.it}}}

\markboth{DRAFT}{Shell \MakeLowercase{\textit{et al.}}: Bare Demo
of IEEEtran.cls for Journals}

\maketitle

\begin{abstract}
In this technical note, we deal with a spectrum approximation problem arising in THREE-like multivariate spectral estimation
approaches. The solution to the problem minimizes a suitable divergence index with respect to an {\em a priori} spectral density.
We derive a new divergence family between multivariate spectral densities which takes root in the prediction theory.
Under mild assumptions on the {\em a priori} spectral density, the approximation problem, based on this new divergence family, admits a family of
solutions. Moreover, an upper bound on the complexity degree of these solutions is provided.
\end{abstract}


\section{Introduction}

The growing interest in the last years on the THREE-like multivariate spectral estimation approach, see the pioneering work \cite{A_NEW_APPROACH_BYRNES_2000}, is certainly due to its appealing features.
First, it leads to a convex optimization problem whose solution is a spectral density typically having an upper bound on the complexity degree.
Second, the numerical solution is computed by efficient algorithms whose global convergence is guaranteed. Third, such approach allows
a wide degree of freedom, through the design of a bank of filters, in order to allow an higher resolution of the solution in the frequency bands of interest.

A THREE-like method may be outlined as follows. A finite length sequence extracted from a realization of a stochastic process, say $y$, is fed to a bank of filters.
The output covariance is then used to extract information on the process by considering the family of spectral densities matching such output covariance.
Therefore, the estimate of the spectral density of $y$ is chosen in this family. This task is accomplished by solving a spectrum approximation problem whose solution minimizes a divergence index with respect to
an {\em a priori} spectral density.

It is clear that the solution highly depends on the divergence
index. Many divergence indexes has been proposed in the last
decade,
\cite{KL_APPROX_GEORGIUO_LINDQUIST,RELATIVE_ENTROPY_GEORGIOU_2006,Hellinger_Ferrante_Pavon,FERRANTE_TIME_AND_SPECTRAL_2012,BETA,ALPHA}.
In particular, in \cite{BETA} a multivariate extension to the Beta
divergence family, \cite{BASU_ROBUST_1998}, has been introduced.
Making additional assumptions on the {\em a priori} spectral
density besides bounded {\em McMillan} degree, it is possible to
show that the corresponding spectrum approximation problem leads
to a family of solutions. Moreover, such a family of solutions
connects the ones obtained with the multivariate {\em
Itakura-Saito} distance, \cite{FERRANTE_TIME_AND_SPECTRAL_2012},
and the multivariate {\em Kullback-Leibler} divergence,
\cite{RELATIVE_ENTROPY_GEORGIOU_2006}. It is worth noting that the
(scalar) {\em Itakura-Saito} distance has been derived from the
maximum likelihood of speech spectral densities,
\cite{ITAKURA_SAITO_1968}, whereas the (scalar) {\em
Kullback-Leibler} divergence has been introduced to measure the
difference between two probability distributions,
\cite{KULLBACK_LEIBLER_1951}. Then, in \cite{ALPHA} the family of
solutions obtained by using the scalar Alpha divergence family,
\cite{AMARI_DIFFERENTIAL_GEOMETRIC}, has been considered. This
family connects the solutions obtained by using the {\em
Kullback-Leibler} divergence
\cite{RELATIVE_ENTROPY_GEORGIOU_2006,KL_APPROX_GEORGIUO_LINDQUIST}.
Here, it is only required that the {\em a priori} spectral density
has bounded {\em McMillan} degree. However, such a result cannot
be extended to the multichannel case.

The main result of this note is the definition of a new multivariate divergence
family which compares two spectral densities in the context of optimal prediction.
This family also includes the multivariate {\em Itakura-Saito} distance.
It is interesting to note this divergence family is not even known in the scalar case and it can be derived from the Alpha and the Beta divergence family.
The corresponding spectrum approximation is then tackled. It turns out that it is possible to characterize a family of solutions to the problem with bounded {\em McMillan} degree
by only making the mild assumption that the {\em a priori} spectral density has bounded {\em McMillan} degree as well. Finally, we face the feasibility issue for the spectrum approximation problem.
Indeed, since the output covariance is estimated from the data, the family of spectral densities matching the (estimated) output covariance is typically empty. Hence, we approximate the output covariance, by solving a state covariance estimation problem, in such a way that the spectrum approximation problem is feasible.

The outline of the technical note is as follows. In Section \ref{section_distance} we introduce the new multivariate divergence family and in Section \ref{section_THREE} the
corresponding spectrum approximation problem. In Section \ref{section_state_cov} we deal with the state covariance estimation problem. In Section \ref{section_simulazioni} a simulation study, showing the features of the family of solutions, is then provided.

Finally, we introduce some convention and notation which will be used throughout the note: integration is on the unit circle with respect to the normalized
{\em Lebesgue} measure. A star denotes transposition plus
conjugation, $\Qc_n\subset \Rs^{n \times n}$ denotes the
$n(n+1)/2$-dimensional real space of $n$-dimensional symmetric
matrices and $\Qc_{n,+}$ denotes the corresponding cone of
positive definite matrices.

\section{A divergence family in the context of optimal
prediction}\label{section_distance}
In this section we introduce a
new divergence family for comparing two spectral densities which
takes root in the prediction theory. Let $y=\Set{y_k;\;k\in\Zs}$
be a stationary, full-rank, purely non-deterministic, zero mean,
$\Rs^m$-valued gaussian process completely characterized by the
spectral density $\Phi(\e{j\vartheta})\in \Stm$, where $\Stm$
denotes the family of coercive, bounded,  $\Cs^{m\times m}$-valued
spectral density functions on the unit circle $\Ts$. When
$\Phi(\e{j\vartheta})$ is not known, we fix an {\em a priori}
spectral density $\Psi(\e{j\vartheta})\in\Stm$ describing $y$. Let
$\hat y_k$ be the least-square linear one-step-ahead predictor
based on $\Psi$, and $e_k:=y_k-\hat y_k$  the corresponding
innovation process. Accordingly $e^N_k:=L_\Omega^{-1}(y_k-\hat
y_k)$ represents the normalized innovation process, with $L_\Omega
L_\Omega^T=\Es[e_ke_k^T]$. Since $\Psi\in\Stm$, it admits a unique
canonical left spectral factor such that \eq
\Psi(\e{j\vartheta})=W_\Psi(\e{j\vartheta})W_\Psi(\e{j\vartheta})^*\eeq
with $W_\Psi(z)\in\mathrm{H}_2^{m\times m}(\Ds)$, $\det
W_{\Phi}(z)\neq 0$ in $\Ds:=\Set{z\;:\; |z|\geq 1}$, and
$W_\Psi(\infty)=L_\Omega$. Here, $\Hs_2^{m\times m}(\Ds)$ denotes
the {\em Hardy} space of analytic functions in $\Ds$ with square
integrable radial limits. Thus, the normalized innovation process
$e_k^N$ is obtained by filtering $y$ through the whitening filter,
say $W_\Psi^{-1}(z)$, \cite{KAILATH_LIN_EST}: \eq
e^N_k=\sum_{l=-\infty }^k a_{\Psi}(k-l)y_{l},\eeq with \eq
W_\Psi^{-1}(z)=\sum_{l=0}^\infty a_{\Psi}(l)z^{-l}, \;\; z\in\Ds.
\eeq The spectral density of the normalized innovation process is
\eq
E^N(\e{j\vartheta})=W_\Psi(\e{j\vartheta})^{-1}\Phi(\e{j\vartheta})W_\Psi(\e{j\vartheta})^{-*}.\eeq
Clearly, if the {\em a priori} spectral density $\Psi$ coincides
with the true spectral density $\Phi$, we have
$E^N(\e{j\vartheta})=I$, i.e. the normalized innovation process is
white gaussian noise (WGN) with zero mean and variance $I$.
Therefore,  $E^N$ represents a mismatch criterium which naturally
occurs in prediction error estimation,
\cite{LINDQUIST_PREDICTION_ERROR_2007,LINDQUIST_PICCI}.

This lead us, as suggested in \cite{GEORGIOU_DISTANCES_PSD_2012},
to measure the mismatch between the true spectral density $\Phi$
and the {\em a priori} spectral density $\Psi$ by quantifying the
mismatch between $E^N$ and $I$: \eq \label{def_dist_generica}
\Sc(\Phi\|\Psi)=\int \Fc(W_\Psi^{-1}\Phi W_\Psi^{-*}),\eeq where
$\Fc:\Qc_{m,+}\rightarrow \Rs$ is a suitable continuous function
such that $\Fc(P)\geq 0$ and equality holds if and only if $P=I$.
Here, we consider the following function parameterized by
$\tau\in\Rs\setminus\{0,1\}$: \eq
\label{def_d_tau}\Fc_\tau(P)=\tr\left[\frac{1}{\tau(\tau-1)}P^\tau-\frac{1}{\tau-1}P\right]+\frac{m}{\tau}.\eeq
Substituting (\ref{def_d_tau}) in (\ref{def_dist_generica}), we
get the following divergence family: \eqn
\label{tau_divergence}\Sc_T^{(\tau)}(\Phi\|\Psi):= \int\tr \left[
\frac{1}{\tau(\tau-1)}\left(W_\Psi^{-1}\Phi
W_\Psi^{-*}\right)^{\tau}\right.\nn\\
\left. -\frac{1}{\tau-1}\Phi
\Psi^{-1}\right]+\frac{m}{\tau}.\eeqn
\rem The multivariate Alpha and Beta divergence family, \cite{ALPHA,BETA}, are defined as follows, respectively:
 \eqn \Sc_A^{(\alpha)}(\Phi\|\Psi)&:=& \int\tr\left[
\frac{1}{\alpha(\alpha-1)}\Phi^{\alpha}\Psi^{1-\alpha}-\frac{1}{\alpha-1}\Phi+\frac{1}{\alpha}\Psi\right]\nn\\
\Sc_B^{({\beta})}(\Phi\|\Psi)&:=& \int\tr
\left[\frac{1}{\beta(\beta-1)}\Phi^{\beta}-\frac{1}{\beta-1}\Phi\Psi^{\beta-1}+\frac{1}{\beta}\Psi^\beta\right]
\nn\eeqn where $\alpha,\beta\in\Rs\setminus\{0,1\}$. There exists
a connection among the above divergences and $\Sc_T^{(\tau)}$: \eq
\Sc_T^{(\tau)}(\Phi\|\Psi)=\Sc_A^{(\tau)}(W_\Psi^{-1}\Phi
W_\Psi^{-*}\| I)=\Sc_B^{(\tau)}(W_\Psi^{-1}\Phi W_\Psi^{-*}\|
I).\eeq Namely, $\Sc_T^{(\tau)}$ measures the dissimilarity
between $E^N$ and $I$ through $\Sc_A^{(\tau)}$ or equivalently
$\Sc_B^{(\tau)}$. \erem
 \prop The following facts
hold:\begin{enumerate}
    \item $\Sc_T^{(\tau)}$ can be extended by continuity for $\tau=0$
    and $\tau=1$:  \eqn &&  \lim_{\tau\rightarrow 0} \Sc_T^{(\tau)}(\Phi\|\Psi)=\Sc_{\mathrm{0}}(\Phi\|\Psi)\nn\\  && \lim_{\tau\rightarrow 1}\Sc_T^{(\tau)}(\Phi\|\Psi)= \Sc_{\mathrm{1}}(\Phi\|\Psi)\eeqn
    where
    \eqn \label{def_IS_KL} && \Sc_{\mathrm{0}}(\Phi\|\Psi):= \int \tr\left[\log\Psi-\log \Phi+\Phi\Psi^{-1}\right] -m\nn\\
    && \Sc_{\mathrm{1}}(\Phi\|\Psi):= \int\tr\left[ W_\Psi^{-1}\Phi
W_\Psi^{-*}\log(W_\Psi^{-1}\Phi W_\Psi^{-*})\right.\nn\\
&& \hspace*{2cm}\left. -\Phi \Psi^{-1}\right] +m.
\eeqn
    \item $\Sc_T^{(\tau)}(\cdot\|\Psi)$ is strictly convex over
    $\Stm$
    \item $\Sc_T^{(\tau)}(\Phi\|\Psi)\geq0$ and equality holds if and only if
    $\Phi=\Psi$.
\end{enumerate}
\eprop Here, $\Sc_{\mathrm{0}}$ is the multivariate {\em
Itakura-Saito} distance, \cite{FERRANTE_TIME_AND_SPECTRAL_2012},
between $\Phi$ and $\Psi$. Recalling that the multivariate {\em
Kullback-Leibler} divergence extended to spectral densities with
different zeroth-moment is defined as, \cite{BETA}, \eq
\Sc_{\mathrm{KL}}(\Omega_1\|\Omega_2)=\int
\tr[\Omega_1(\log(\Omega_1)-\log(\Omega_2))-\Omega_1+\Omega_2]\eeq
we conclude that $\Sc_{\mathrm{1}}$ is the multivariate {\em
Kullback-Leibler} divergence between $E^N$ and $I$.

{\em Sketch of the Proof:} Since $\Psi\in\Stm$, the linear map
$f_\Psi: \Phi\mapsto W_\Psi^{-1}\Phi W_\Psi^{-*}$ is bijective.
Moreover, $\Sc_T^{(\tau)}(\Phi\|\Psi)=\Sc_B^{(\tau)}(\cdot \|
I)\circ f_\Psi(\Phi)$. Thus, it is sufficient to apply Proposition
3.2 and Proposition 3.3 in \cite{BETA} to get the statement. \qed

$\Sc_T^{(\tau)}$ has been defined in (\ref{tau_divergence})
through the canonical left spectral factor $W_\Psi$ of $\Psi$, however it
does not depend on this particular choice. \prop Let
$\overline{W}_\Psi$ be any left square spectral factor of $\Psi$.
Define\eqn && \overline{\Sc}_T^{(\tau)}(\Phi\|\Psi)=
 \int\tr\left[\frac{1}{\tau(\tau-1)} (\overline{W}_\Psi^{-1} \Phi\overline{W}_\Psi^{-*})^\tau\right. \nn\\ && \hspace{3.8cm}\left. -\frac{1}{\tau-1}\Phi \Psi^{-1}\right]+\frac{m}{\tau}\eeqn where $\tau\in\Rs\setminus \{0,1\}$. Then, \eq \label{equivalence_relation}\Sc_T^{(\tau)}(\Phi\|\Psi)=\overline{\Sc}_T^{(\tau)}(\Phi\|\Psi)\eeq
and this equality also holds for $\tau\rightarrow 0$ and $\tau\rightarrow 1$.\eprop

\IEEEproof The spectral factor $\overline{W}_\Psi$ can be obtained through the canonical one as follows \eq \overline{W}_\Psi= W_\Psi U^*\eeq
where $U$ is an $m\times m$ all pass function, i.e. $UU^*=I$ on $\Ts$. Let $\overline{W}_\Psi^{-1} \Phi\overline{W}_\Psi^{-*}=VDV^*$ be a pointwise SVD of $\overline{W}_\Psi^{-1} \Phi\overline{W}_\Psi^{-*}$, therefore $D$ is diagonal and $VV^*=I$ on $\Ts$. Hence,
\eqn  \overline{W}_\Psi^{-1} \Phi\overline{W}_\Psi^{-*}&=&(W_\Psi U^*)^{-1} \Phi (W_\Psi U^*)^{-*}\nn\\
&=& U W_\Psi^{-1} \Phi W_\Psi^{-*} U^*=UV D V^* U^* \eeqn
and $UV$ is an $m \times m$ all pass function, in fact \eq (UV)(UV)^*=UVV^*U^*=UU^*=I.\eeq
Thus, $VDV^*$ and $(UV)D(UV)^*$ are two pointwise SVD of $\overline{W}_\Psi^{-1} \Phi\overline{W}_\Psi^{-*}$.
In order to show that (\ref{equivalence_relation}) holds for $\tau\in\Rs\setminus\{0,1\}$, it is sufficient to show that
$\tr[(W_{\Psi}^{-1}\Phi W_{\Psi}^{-*})^\tau]=\tr[(\overline{W}_{\Psi}^{-1}\Phi \overline{W}_{\Psi}^{-*})^\tau]$:
\eqn && \tr\left[ (\overline{W}_\Psi^{-1} \Psi\overline{W}_\Psi^{-*})^\tau\right]=\tr\left[ (UV D V^*U^*)^\tau\right]\nn\\ &&\hspace{0.5cm}=\tr\left[ UV D^\tau V^*U^*\right]=\tr\left[ V D^\tau V^*\right]\nn\\ && \hspace{0.5cm}=\tr\left[ (V D V^*)^\tau \right]= \tr\left[(W_\Psi^{-1} \Psi W_\Psi^{-*})^\tau\right] \eeqn
accordingly
$\Sc_T^{(\tau)}(\Phi\|\Psi)=\overline{\Sc}_T^{(\tau)}(\Phi\|\Psi)$.
For the limit cases, we have
\eqn && \lim_{\tau\rightarrow 0} \Sc_{\tau}(\Phi\|\Psi)=\int\tr\left[\log\Psi-\log\Phi+\Phi\Psi^{-1}\right]-m.\nn\\
&& \lim_{\tau\rightarrow 1}\Sc_{\tau}(\Phi\|\Psi)= \int\tr\left[\overline{W}_\Psi^{-1} \Psi\overline{W}_\Psi^{-*}\log( \overline{W}_\Psi^{-1} \Psi\overline{W}_\Psi^{-*})\right.\nn\\ && \hspace{5.0cm}\left.-\Phi\Psi^{-1}\right]+m.\eeqn
Therefore, (\ref{equivalence_relation}) holds for $\tau\rightarrow 0$. For $\tau\rightarrow 1$ we have \eqn && \tr\left[ \overline{W}_\Psi^{-1} \Psi \overline{W}_\Psi^{-*}\log(\overline{W}_\Psi^{-1} \Psi \overline{W}_\Psi^{-*})\right]\nn\\ && \hspace{0.5cm}=\tr\left[UVDV^*U^*UV\log(D)V^*U^* \right]\nn\\
&& \hspace{0.5cm}=\tr\left[VDV^*V\log(D)V^* \right]\nn\\ &&\hspace{0.5cm}= \tr\left[W_\Psi^{-1} \Psi W_\Psi^{-*}\log(W_\Psi^{-1} \Psi W_\Psi^{-*})\right]\eeqn
therefore (\ref{equivalence_relation}) holds  for $\tau\rightarrow 1$.  \qed

It is worth noting that $\Sc_T^{(\tau)}$ remains invariant under
congruence, that is \eq  \Sc_T^{(\tau)}(\Phi\|\Psi)=
\Sc_T^{(\tau)}(T^*\Phi T\| T^*\Psi T)\eeq for any $\Cs^{m\times
m}$-valued function $T(z)$ invertible on $\Ts$. This is a natural
property to demand since it implies that the divergence between
spectral densities does not change under coordinate
transformation. In particular, $\Sc_T^{(\tau)}$ is invariant to
scaling: such a property seems essential, for  instance, for
speech and image systems, due to an apparent agreement with
subjective qualities of sound and images. This invariance
property, however, does not hold for $\Sc_A^{(\alpha)}$ and
$\Sc_B^{(\beta)}$. Finally, it is worth noting that
$\Sc_T^{(\tau)}$ is not even known in the scalar case (i.e.
$m=1$).

\section{THREE-like spectral estimation}\label{section_THREE}

Consider the stochastic process $y$ of Section
\ref{section_distance} with unknown spectral density
$\Phi(\e{j\vartheta})\in\Stm$. Assume that, the given {\em a
priori} spectral density $\Psi\in\Stm$ has bounded {\em McMillan}
degree. Then, a finite length sequence $\mathbf{y}_1 \ldots
\mathbf{y}_N$ extracted from a realization of $y$ is available. We
want to find an estimate of $\Phi$ by exploiting $\Psi$ and
$\mathbf{y}_1 \ldots \mathbf{y}_N$. According to the THREE-like
approach, we design a rational filter \eq \label{filter_bank}
G^\prime(z)=(zI-A)^{-1}B,\eeq where $A\in\Rs^{n\times n}$ is a
stability matrix, $B\in\Rs^{n\times m}$ is full rank with $n> m$
and $(A,B)$ is a reachable pair. We compute an estimate
$\hat{\Sigma}\succ 0$, based on the data $\mathbf{y}_1 \ldots
\mathbf{y}_N$, of the steady state covariance $\Sigma=\Sigma^T
\succ0$ of the state $x_k$ of the filter \eq
x_{k+1}=Ax_k+By_k.\eeq Then, an estimate of $\Phi$ is given by
solving the following spectrum approximation problem. \pb
\label{pb_spectrum_approx}Given $\Psi\in\Stm$,
$G(z)=\hat{\Sigma}^{-\frac{1}{2}}(zI-A)^{-1}B$, and
$\nu\in\Ns_+:=\Ns\setminus \{0\}$,\eqn \label{set_approx_pb}&&
\mathrm{minimize}\; \;\Sc_T^{(1-\nu^{-1})}(\Phi\|\Psi)\;\;
\mathrm{over \;the\; set} \nn\\&& \hspace{1cm} \; \;
\Jc=\Set{\Phi\in\Stm\;|\;\int G\Phi G^*=I }.\eeqn \epb  The
constraint in (\ref{set_approx_pb}) is equivalent to \eq \int
G^\prime \Phi (G^\prime)^*=\hat\Sigma.\eeq Accordingly, the
optimal solution to Problem \ref{pb_spectrum_approx} (if it does
exist) matches the data encoded by $\hat\Sigma$ and is such that
the one-step-ahead predictor based on $\Psi$ is as close as
possible to be the optimal one. Note that, the set $\Jc$ depends
on the estimate $\hat \Sigma$. In this section we assume that
$\hat \Sigma$ is chosen in such a way that $\Jc$ is non-empty,
that is Problem \ref{pb_spectrum_approx} is feasible. In Section
\ref{section_state_cov}, we will show how to compute such a
$\hat\Sigma$.

In \cite{FERRANTE_TIME_AND_SPECTRAL_2012}, it was already shown
that Problem \ref{pb_spectrum_approx} admits a unique solution for
$\nu=1$. Thus, we deal with the case $\nu>1$. Since Problem
\ref{pb_spectrum_approx} is a constrained convex optimization
problem, we consider the corresponding {\em Lagrange}
functional\eqn && \hspace{-0.9cm} L_\nu(\Phi, \Lambda)\nn\\ &&
\hspace{-0.5cm}=\Sc_T^{(1-\nu^{-1})}(\Phi\|\Psi) -m
\frac{\nu}{\nu-1}+\Sp{\int G\Phi G^*-I}{\Lambda}\nn\\
&& \hspace{-0.5cm}= \int \tr\left[
\frac{\nu^2}{1-\nu}(W_\Psi^{-1}\Phi
W_\Psi^{-*})^{\frac{\nu-1}{\nu}}+\nu (W_\Psi^{-1}\Phi
W_\Psi^{-*})\right.\nn\\ && \left. +G^*\Lambda
G\Phi\right]-\tr[\Lambda]\eeqn where we exploited the fact that
the term $m\frac{\nu}{\nu-1}$ plays no role in the optimization
problem. Note that, the domain of $L_\nu(\cdot,\Lambda)$ is
$\Sc_+^m(\Ts)$ and $\Lambda\in\Qc_n$ is the {\em Lagrange
multiplier} associated to the constraint in (\ref{set_approx_pb}).
Consider the vector space\eq \Qc^G_n:=\left\{ \int G \Phi G^*
\hbox{ s.t. } \Phi\in \Vc(\Stm)\right\}\eeq where $\Vc(\Stm)$
denotes the vector space generated by $\Stm$. In
\cite{MATRICIAL_ALGTHM_RAMPONI_FERRANTE_PAVON_2009}, it was shown
that $\Lambda$ can be uniquely decomposed as
$\Lambda_G+\Lambda_\bot$ where $\Lambda_G\in\Qc_n^G$ and
$\Lambda_\bot\in(\Qc_n^G)^{\bot}$. Moreover, $G^*(\e{i\vartheta})
\Lambda_\bot G(\e{i\vartheta})\equiv 0$ and $\tr[\Lambda_\bot]=0$.
Accordingly, $\Lambda_\bot$ does not affect the {\em Lagrange}
functional and we can restrict $\Lambda\in\Qc_n^G$. Since
$L_\nu(\cdot,\Lambda)$ is strictly convex over $\Sc_+^m(\Ts)$, it
is sufficient to show the existence of a stationary point for
$L_\nu(\cdot,\Lambda)$ in order to prove the existence of a
(unique) solution to the unconstrained minimization problem
$\Phi_\nu(\Lambda):=\arg\Min{\Phi}\Set{L_\nu(\Phi,\Lambda)\;|\;
\Phi\in\Sc_+^m(\Ts)}$. The first variation of
$L_\nu(\cdot,\Lambda)$ in each direction $\delta\Phi\in L^{m
\times m}_\infty(\Ts)$ is: \eqn
\label{variazione_delta_phi_beta}&&\hspace{-0.5cm} \delta
L_\nu(\Phi,\Lambda;\delta \Phi)= \int\tr
\left[-\nu(W_\Psi^{-1}\Phi
W_\Psi^{-*})^{-\frac{1}{\nu}}W_\Psi^{-1}\delta\Phi
W_\Psi^{-*}\right. \nn\\ &&\hspace{0.2cm}\left.+\nu
W_\Psi^{-1}\delta\Phi W_\Psi^{-*}+G^*\Lambda G\delta\Phi \right]
\nn\\&& \hspace{-0.2cm}= \int\tr \left[\left(-\nu
W_\Psi^{-*}(W_\Psi^{-1}\Phi
W_\Psi^{-*})^{-\frac{1}{\nu}}W_\Psi^{-1}\right.\right. \nn\\ &&
\hspace{0.2cm}\left.\left.+\nu W_\Psi^{-*}W_\Psi^{-1}+G^*\Lambda
G\right)\delta\Phi \right] \eeqn where we used the expression of
the first variation of the exponentiation of $X\in\Qc_{n,+}$ to
$c\in\Rs$ given in \cite{BETA}:\eq
\label{variation_exp}\delta(\tr[X^c];\delta
X)=c\tr[X^{c-1}],\;\;\delta X\in\Qc_{n} .\eeq Note that, $-\nu
W_\Psi^{-*}(W_\Psi^{-1}\Phi
W_\Psi^{-*})^{-\frac{1}{\nu}}W_\Psi^{-1}+\nu
W_\Psi^{-*}W_\Psi^{-1}+G^*\Lambda G\in L^{m \times
m}_\infty(\Ts)$. Thus, (\ref{variazione_delta_phi_beta}) is zero
$\forall \delta\Phi\in L^{m \times m}_\infty(\Ts)$ if and only if
\eq (W_\Psi^{-1}\Phi W_\Psi^{-*})^{-\frac{1}{\nu}}=I
+\frac{1}{\nu} W_\Psi^{*}G^*\Lambda G  W_\Psi. \eeq Since
$(W_\Psi^{-*}\Phi W_\Psi^{-*})^{-\frac{1}{\nu}}\in\Sc_+^m(\Ts)$,
the set of the admissible {\em Lagrange} multipliers is \eq
\Lp:=\Set{\Lambda\in \Qc^G_{n}\;|\;
I+\frac{1}{\nu}W_\Psi^*G^*\Lambda GW_\Psi\succ 0\hbox{ on
}\Ts}.\eeq Therefore, if $\Lambda\in\Lc_+$ then the unique
stationary point for $L_\nu(\cdot,\Lambda)$ is
\eq\label{ottimo_tau}\Phi_\nu(\Lambda):=W_\Psi\left(I+\frac{1}{\nu}W_\Psi^*G^*\Lambda
GW_\Psi\right)^{-\nu}W_\Psi^*\eeq which coincides with the unique
minimum point for $L_\nu(\cdot,\Lambda)$. \prop
\label{prop_optimal_form} If $\Phi_\nu$ is a minimizer of Problem
\ref{pb_spectrum_approx}, then it has bounded {\em McMillan}
degree which is less than or equal to $\nu(\deg[\Psi]+2 n)$.
Moreover the following facts hold:
\begin{enumerate}
  \item Among all the spectral densities $\Phi_\nu$ with $\nu\in\Ns_+$, the spectral density with the smallest upper bound on the {\em McMillan} degree corresponds to
  $\Sc_0$ defined in (\ref{def_IS_KL})
  \item As $\nu\rightarrow + \infty$, $\Phi_\nu$ converges to the optimal form corresponding to $\Sc_1$ defined in (\ref{def_IS_KL}).
\end{enumerate} \eprop
\IEEEproof In \cite{FERRANTE_TIME_AND_SPECTRAL_2012}, it was shown that (\ref{ottimo_tau}) holds for $\nu=1$ and $\deg[\Phi_1]\leq \deg[\Psi]+2n$. For the case $\nu>1$, $\Phi_\nu(\Lambda)= L_\Lambda
Q_\Lambda^{\nu-2}L_\Lambda^T$ where
$Q_\Lambda:=\left(I+\frac{1}{\nu}W_\Psi^*G^*\Lambda G W_\Psi
\right)^{-1}$, $L_\Lambda:=(W_\Psi^{-1}+W_\Psi^*G^*\Lambda
G)^{-1}$, $G$ has bounded {\em McMillan} degree by assumption, and
$W_\Psi$ has bounded degree because is the canonical left spectral
factor of a spectral density with bounded degree. Since
$\deg[Q_\Lambda]$ and $\deg[L_\Lambda]$ are less than or equal to
$\deg[\Psi]+2n$, we conclude that $\Phi_\nu(\Lambda)$ has
bounded {\em McMillan} degree which is less than or equal to $\nu(\deg[\Psi]+2 n)$.\\
{\em Point 1.} Since $\nu(\deg[\Psi]+2 n)$ is an increasing function in $\nu\in\Ns_+$, its minimum is achieved with $\nu=1$, i.e.
with $\Sc_T^{(0)}=\Sc_{\mathrm{0}}$.\\
{\em Point 2.} It is not difficult to show that the optimal form
obtained by using $\Sc_{\mathrm{1}}$ in Problem
\ref{pb_spectrum_approx} is $\Phi_\infty(\Lambda):=W_\Psi
\e{-W^*_\Psi G^*\Lambda GW_\Psi}W^*_\Psi$. By using the limit, see
\cite[Proposition 4.1]{BETA}: \eq
\lim_{\nu\rightarrow\infty}\left(I+\frac{1}{\nu}
X\right)^{-\nu}=\e{-X},\;\; X\in\Qc_n, \eeq we get \eq
\label{limit_opt_form}\lim_{\nu\rightarrow\infty}
\Phi_\nu(\Lambda)=W_\Psi \e{-W_\Psi^*G^*\Lambda G W_\Psi}
W_\Psi^*=\Phi_\infty(\Lambda).\eeq  \qed \rem In \cite{BETA}, the
multivariate  Beta divergence family $\Sc_B^{(1-\nu^{-1})}$ with
$\nu\in\Ns_+$ has been considered. The optimal form of the
spectrum approximation problem is \eq
\Phi_{B,\nu}(\Lambda):=\left(\Psi^{-\frac{1}{\nu}}+\frac{1}{\nu}G^*\Lambda
G\right)^{-\nu}\eeq and the assumption that $\Psi$ has bounded
{\em McMillan} degree is not sufficient to guarantee that
$\Psi^{\frac{1}{\nu}}$, and also $\Phi_{B,\nu}$, has bounded
degree. In \cite{ALPHA}, the scalar Alpha divergence family
$\Sc_A^{(1-\nu^{-1})}$ with $\nu\in\Ns_+$ has been considered. The
corresponding optimal form is \eq
\Phi_{A,\nu}(\Lambda):=\frac{\Psi}{(I+\frac{1}{\nu}G^*\Lambda
G)^\nu}\eeq and $\Psi$ with bounded {\em McMIllan} degree implies
that $\Phi_{A,\nu}$ has bounded degree. However, such a result
cannot be extended to the multichannel case, i.e. $m>1$. In view
of Proposition \ref{prop_optimal_form}, $\Sc_T^{(1-\nu^{-1})}$
with $\nu\in\Ns_+$ is the unique divergence family always leading
to a bounded degree optimal form in the multichannel case once
$\Psi$ has bounded degree.\erem

Since $\Phi_\nu$ is the unique minimum
point for $L_\nu(\cdot,\Lambda)$ over $\Stm$, if we produce
$\Lambda^\circ\in\Lp$ such that $\int
G\Phi_\nu(\Lambda^\circ)G^*=I$, then $\Phi_\nu(\Lambda^\circ)$ is
the unique solution to Problem \ref{pb_spectrum_approx}. To this end, consider the dual functional
defined over $\Lp$: \eq J_\nu(\Lambda)=  \frac{\nu}{1-\nu} \int\tr\left[
\left(I+\frac{1}{\nu}W^*_\Psi G^*\Lambda
GW_\Psi\right)^{1-\nu}\right]-\tr[\Lambda].\eeq
 \teo \label{prop_dual} If $\Jc$ is a non-empty set, then the dual problem $\max\Set{J_{\nu}(\Lambda)\;|\; \Lambda
 \in\Lp}$ with $\nu>1$ admits a unique solution. Such a solution, say $\Lambda^\circ$, satisfies \eq \label{constraint_lambda_circ} \int
G\Phi_\nu(\Lambda^\circ)G^*=I.\eeq\eteo \IEEEproof Recall that $Q_\Lambda=(I+\frac{1}{\nu} W_\Psi^*G^*\Lambda G W_\Psi)^{-1}$. By similar argumentations used in \cite[Theorem 5.1]{BETA}, it is possible to show that $J_\nu\in\Cc^2(\Lc_+)$, and
\eqn && \hspace{-0.6cm}\delta J_\nu(\Lambda;\delta\Lambda)\nn\\ &&\hspace{-0.4cm} =\tr \int\left[ Q_\Lambda^{\nu} W_\Psi^* G^* \delta \Lambda G W_\Psi \right]-\tr[\delta \Lambda]  \label{prima_var_J}\\
&& \hspace{-0.6cm}\delta^2 J_\nu(\Lambda;\delta \Lambda)\nn\\ && \label{seconda_var_J}\hspace{-0.4cm}= -\frac{1}{\nu}\sum_{l=1}^{\nu} \int \tr\left[Q_\Lambda^l W_\Psi^* G^* \delta \Lambda G  W_\Psi   Q_\Lambda^{\nu+1-l}\right. \nn\\ && \hspace{-0.1cm} \left. \times W_\Psi^* G^*\delta \Lambda G W_\Psi   \right]\eeqn
where we used the fact that the first variation of $Q_\Lambda$ in direction $\delta \Lambda\in\Qc_n$ is
\eq Q_{\Lambda;\delta \Lambda}=-\frac{1}{\nu}Q_\Lambda W_\Psi^* G^* \delta\Lambda G W_\Psi Q_\Lambda\eeq
and the first variation of the map $\Ic: A\mapsto A^\nu$ in direction $\delta A\in\Qc_n$ is: \eq \delta (\Ic(A); \delta A)=\sum_{l=1}^\nu A^{l-1} \delta A A^{\nu-l}.\eeq
Since $\nu>1$ and the trace of the integrands in (\ref{seconda_var_J}) is nonnegative, we have $\delta^2 J_\nu\leq 0$. If $\delta^2 J_\nu=0$, then $G^*\delta \Lambda G\equiv 0$, namely $\delta \Lambda\in (\Qc_n^G)^\bot$, see \cite{MATRICIAL_ALGTHM_RAMPONI_FERRANTE_PAVON_2009}. Since $\delta\Lambda \in\Qc_n^G$, we conclude that $\delta \Lambda=0$.
This means that $\delta^2 J_\nu$ is negative definite. Thus, $J_\nu$ is strictly concave on $\Lc_+$ and the dual problem admits at most one solution $\Lambda^\circ$ which must annihilate (\ref{prima_var_J}) for each $\delta\Lambda$. This implies that $\Lambda^\circ$ satisfies (\ref{constraint_lambda_circ}).

It remains to be shown that $J_\nu$ takes a maximum value on $\Lc_+$. Note that $J_\nu(0)=m\frac{\nu}{1-\nu}$, accordingly we can restrict the search of a maximum point to the nonempty set $\Kc:=\{\Lambda \in\Lc_+ \hbox{ s.t. } J_\nu(\Lambda)\geq J_\nu(0)\}$. By similar argumentations used in \cite[Theorem 5.2]{BETA}, it is possible to show that
\eqn \lim_{\Lambda\rightarrow \partial \Lc_+} J_\nu(\Lambda)=-\infty\nn\\
\lim_{\|\Lambda \|\rightarrow \infty } J_\nu(\Lambda)=-\infty\eeqn
that is $\Kc$ is compact. Finally, since $J_\nu\in\Cc^2(\Lc_+)$, the existence of the solution follows from the {\em Weierstrass}' Theorem.\qed

The optimal solution $\Lambda^\circ$ can be efficiently computed
by using the matricial {\em Newton} algorithm with backtracking
presented in \cite{MATRICIAL_ALGTHM_RAMPONI_FERRANTE_PAVON_2009}.
Here, the {\em Newton step} $\Delta_{\Lambda_i}$ at the $i$-th
iteration with starting point $\Lambda_i$ is given by solving the
linear equation \eqn\label{eq_Newton_step} \frac{1}{\nu}
\sum_{l=1}^{\nu}\int G W_\Psi Q_{\Lambda_i}^l W_\Psi^* G^*
\Delta_{\Lambda_i} G W_\Psi Q_{\Lambda_i}^{\nu+1-l}W_\Psi^*
G^*\nn\\ \hspace{1cm}=\int  G W_\Psi Q_{\Lambda_i}^{\nu} W_\Psi^*
G^*- I.\eeqn
 Concerning its computation, in \cite{MATRICIAL_ALGTHM_RAMPONI_FERRANTE_PAVON_2009} a sensible and efficient method based on spectral factorization techniques has been presented. Finally, it is possible to prove that this algorithm globally converges, in particular the rate of convergence is quadratic during the last stage.

\section{State Covariance Estimation}\label{section_state_cov}

In this section we face the problem of computing $\hat\Sigma$ in
such a way that the set $\Jc$ is non-empty. In
\cite{ME_ENHANCEMENT_FERRANTE_2012}, (see also
\cite{GEORGIUO_THE_STRUCUTRE_2002,Hellinger_Ferrante_Pavon}) it
was shown that $\Jc$ is non-empty if and only if $\hat \Sigma$ is
positive definite and it belongs to the kernel of the linear
operator \eqn V &:& \Qc_n\rightarrow \Qc_n\nn\\ && Q\mapsto
\Pi_B^\bot(Q-AQA^T) \Pi_B^\bot,\eeqn where
$\Pi_B^\bot=I-B(B^TB)^{-1}B^T$. Thus, we can feed the bank of filters
$G^\prime(z)$ with the finite length sequence $\mathbf{y}_1\ldots \mathbf{y}_N$,
obtaining the output data $\mathbf{x}_1\ldots \mathbf{x}_N$. Then an estimate of
$\Sigma$ is given by the sample covariance $\hat
\Sigma_C:=\sum_{k=1}^N \mathbf{x}_k\mathbf{x}_k^T$ which is normally positive
definite, but it may happen $V(\hat \Sigma_C)\neq 0$ especially
when $N$ is not large. Following the same approach presented in
\cite{ME_ENHANCEMENT_FERRANTE_2012,BETA}, an estimate of $\Sigma$
leading to a non-empty set $\Jc$ is given by  finding a new
positive definite estimate $\hat\Sigma$ such that
$V(\hat\Sigma)=0$ and ``close'' to the estimate $\hat\Sigma_C$.

For measuring the closeness between $\hat\Sigma$ and $\hat
\Sigma_C$, we consider the following divergence index between
$P,Q\in\Qc_{n,+}$ with $\tau\in\Rs\setminus\{0,1\}$: \eq
\label{tau_matrix_divergence}\Dc_T^{(\tau)}(P\|Q):=\tr\left[
\frac{1}{\tau(\tau-1)}(L_Q^{-1}PL_Q^{-T})^{\tau}-\frac{1}{\tau-1}(PQ^{-1})\right]+\frac{m}{\tau}.\eeq
Here, $L_Q$ is the {\em Cholesky} decomposition of $Q$, i.e.
$Q=L_QL_Q^T$. Note that $\Dc_T^{(\tau)}$ is a special case of
$\Sc_T^{(\tau)}$: it is sufficient to pick $\Phi(\e{i\vartheta
})=P$ and $\Psi(\e{i\vartheta })=Q$ in order to obtain
(\ref{tau_matrix_divergence}). Accordingly, $\Dc_T^{(\tau)}$ is
strictly convex with respect to the first argument and it can be
extended by continuity for $\tau=0$ and $\tau=1$: \eq
\lim_{\tau\rightarrow 0}
\Dc_T^{(\tau)}(P\|Q)=\Dc_{\mathrm{0}}(P\|Q),\;\;
\lim_{\tau\rightarrow 1}\Dc_T^{(\tau)}(P\|Q)=
\Dc_{\mathrm{1}}(P\|Q)\eeq
 where
    \eqn && \Dc_{\mathrm{0}}(P\|Q):=\left[ \tr\int \log Q-\log P+PQ^{-1}\right] -m\nn\\
    && \Dc_{\mathrm{1}}(P\|Q):=\left[ \tr\int (L_Q^{-1}P
L_Q^{-T})\log(L_Q^{-1}P L_Q^{-T})\right.\nn\\
&& \hspace*{2cm}\left. -L_Q^{-1}P L_Q^{-T}\right] +m. \eeqn Here,
$\Dc_{\mathrm{0}}$ is the {\em Burg matrix divergence} and
$\Dc_{\mathrm{1}}$ is the {\em Umegaki-von Neumann's relative
entropy} (extended to non-equal trace matrices) between $L_Q^{-1}
P L_Q^{-T}$ and $I$.

We are now ready to introduce the optimization problem for finding $\hat\Sigma$.
\pb \label{pb_approx_covarianze} Given $\hat\Sigma_C\succ 0$ and $\nu\in\Ns_+$, solve \eqn \hbox{minimize } \Dc_T^{(1-\nu^{-1})} (\hat\Sigma \|\hat\Sigma_C) \hbox{ over the set }\nn\\
\hspace{2cm} \Set{\hat\Sigma\in\Qc_{n,+}\;|\; V(\hat \Sigma)=0}.\eeqn\epb

\teo Problem \ref{pb_approx_covarianze} admits a unique solution.\eteo
\IEEEproof  Since $\Dc_T^{(1-\nu^{-1})}(\cdot\| \hat\Sigma_C)$ is strictly convex on $\Qc_{n,+}$, Problem \ref{pb_approx_covarianze} admits at most one solution. Then, the existence of such a solution can be proved by duality theory as in \cite[Section VIII]{BETA}.\qed\\

Also in this case, a globally convergent matricial Newton
algorithm for computing $\hat\Sigma$ may be used, see
\cite{ME_ENHANCEMENT_FERRANTE_2012}. Finally, it is worth noting
that it is also possible to estimate $\Sigma$ through the
estimates of the covariance lags of the input process $y$,
\cite{ON_THE_ESTIMATION_ZORZI_2012,FERRANTE_STRUCTURED_2012}.

\section{Simulation study}\label{section_simulazioni}

In this section we want to test the features of the family of estimators $\Phi_\nu$ through a simulation study.
Assume that a finite length sequence $\mathbf{y}_1\ldots \mathbf{y}_{N}$ is extracted from a realization of a bivariate process $y$ with spectral density $\Phi\in\Stm$. We want to compute the estimates $\Phi_\nu$ with $\nu=1$ and $\nu=2$ of $\Phi$. In view of the results of
Section \ref{section_THREE} and Section \ref{section_state_cov}, we consider the following identification procedure:
\begin{itemize}
\item Choose a low order {\em a priori} spectral density $\Psi$ with bounded {\em McMillan} degree;
  \item Choose $G^\prime(z)$ as in Section \ref{section_THREE}, and $\nu\in\Ns_+$;
  \item Feed the bank of filters $G^\prime$ with the sequence $\mathbf{y}_1\ldots \mathbf{y}_N$, collect the output data $\mathbf{x}_1\ldots \mathbf{x}_N$ and compute $\hat\Sigma_C=\frac{1}{N}\sum_{k=1}^N \mathbf{x}_k\mathbf{x}_k^T$;
  \item Compute $\hat\Sigma_\nu$ by solving Problem \ref{pb_approx_covarianze};
  \item Compute $\Phi_\nu$ by solving Problem \ref{pb_spectrum_approx} with the chosen $\Psi$ and $G(z)=\hat\Sigma_\nu^{-\frac{1}{2}}G^\prime(z)$.
\end{itemize} In the above procedure we have two degrees of freedom: the way for choosing $\Psi$ and the structure of
$G^\prime(z)$.

The {\em a priori} spectral density $\Psi$ may be derived from given laws (e.g. physical laws if $y$ represents a physical phenomenon)
or from the data by applying a (simple) identification method. In both cases $\Psi$ represents a coarse, low order, estimate of $\Phi$.
In our case, $\Psi$ is an ARMA(1,1) which is computed from $\mathbf{y}_1\ldots \mathbf{y}_{N}$ by applying the \textrm{MATLAB}'s  \textrm{PEM} identification method.

Concerning the design of $G^\prime(z)=(zI-A)^{-1}B$, a higher resolution can be
attained by selecting poles in the proximity of the unit circle,
with arguments in the range of frequency of interest,
\cite{A_NEW_APPROACH_BYRNES_2000}. Here, we choose \eq A=\left[
                                                                       \begin{array}{cccc}
                                                                         0 & I & 0 & 0   \\
                                                                         0 & 0 & I & 0 \\
                                                                         0 & 0 & 0 & I  \\
                                                                         0 & 0 & 0 & 0   \\
                                                                       \end{array}
                                                                     \right] \;\; B=\left[
                                                                                      \begin{array}{c}
                                                                                        0 \\
                                                                                        0 \\
                                                                                        0 \\
                                                                                        I \\
                                                                                      \end{array}
                                                                                    \right]\eeq
which implies that $\Sigma$ is a {\em block Toeplitz} matrix whose block diagonals contain the lags $\mathbb{E}[y_k y_{k+j}^T]$ $j=0,1,2,3$ of the covariance function of $y$. Accordingly, Problem \ref{pb_spectrum_approx} becomes a covariance extension problem.

Finally, the sequence $\mathbf{y}_1 \ldots \mathbf{y}_N$ is generated by feeding a bivariate WGN with zero mean and variance $I$ to a (stable) square shaping filter of order 40. The latter is constructed with random coefficients.

\begin{figure}[htbp]
\begin{center}
\includegraphics[width=\columnwidth]{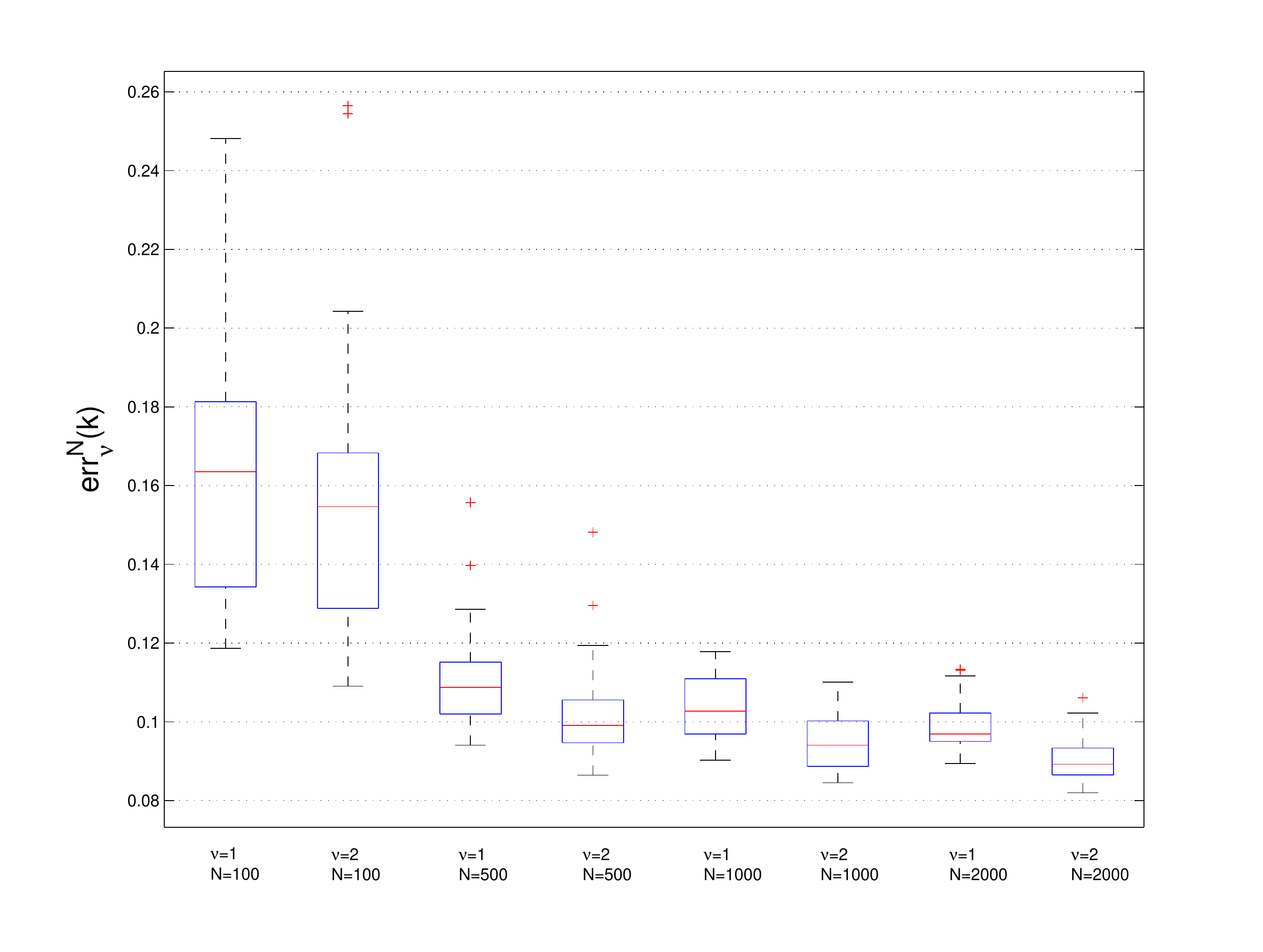}
\end{center}
 \caption{Boxplots of the relative errors averaged over the unit circle for different values of $N$, and $\nu=1$ and $\nu=2$.}\label{figura_boxplots}
\end{figure}

\begin{figure}[htbp]
\begin{center}
\includegraphics[width=\columnwidth]{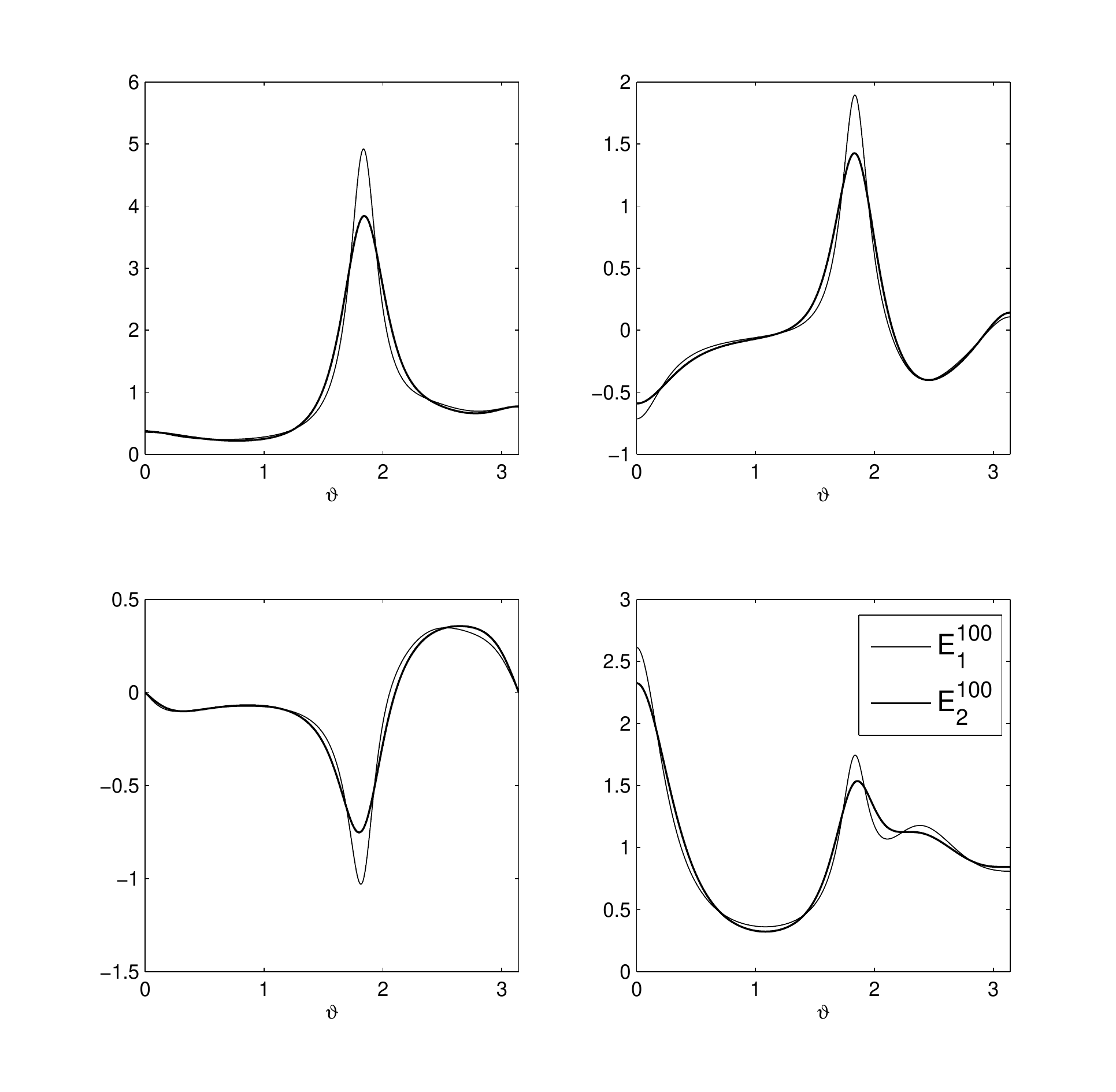}
\end{center}
 \caption{Spectral densities of the normalized innovation process averaged over 50 runs with $N=100$.}\label{figura_innovation_pro}
\end{figure}

We consider four different lengths for the sequence: $N=100$,
$N=500$, $N=1000$ and $N=2000$, i.e. we start from short
observation records up to long observation records. For each $N$,
in order to obtain a comparison reasonably independent of the
specific data set, we perform 50 independent runs where the
sequence $\mathbf{y}_1 \ldots \mathbf{y}_N$ changes in each run.
In this way, we obtain 50 different estimates $\Phi_{\nu,k}^N$
$k=1\ldots 50$. We define then \eq \mathrm{err}_\nu^N(k):=\int
\frac{\| \Phi_{\nu,k}^N-\Phi\|}{\|\Phi\|} \eeq where $\|\cdot \|$
denotes the spectral norm. This is understood as the relative
estimation error of $\Phi_{\nu,k}^N$ averaged over the unit
circle. In Figure \ref{figura_boxplots}  boxplots of the averaged error achieved by
the estimates  with $\nu=1$ and $\nu=2$ for different values of
$N$ is depicted. Clearly, the larger $N$ is the better estimates
are. Moreover, the estimator with $\nu=2$ always outperforms the
one with $\nu=1$. On the other hand, the {\em McMillan} degree of
$\Phi_\nu$ increases by increasing $\nu$: the spectral factor
(i.e. the model) with $\nu=1$ has degree 9, whereas with $\nu=2$
has degree $18$. The spectral density of the normalized innovation
process corresponding to $\Phi_{\nu,k}^N$ is
$E_{\nu,k}^N:=W_\Psi^{-1}\Phi_{\nu,k}^NW_\Psi^{-*}$. We define \eq
E_\nu^N=\frac{1}{50} \sum_{k=1}^{50} E^N_{\nu,k}\eeq which is the
spectral density of the normalized innovation process averaged
over 50 runs and $N$ is fixed. In Figure
\ref{figura_innovation_pro},  $E_\nu^{100}$ with $\nu=1$ and
$\nu=2$ is depicted. We observe that $E_2^{100}$ is more similar
to WGN with variance $I$ than $E_1^{100}$. It is also interesting to analyze the shape of $E_{\nu,k}^N$ changing $\nu$. For simplicity consider the scalar case, $m=1$, and let
$\bar{e}=\Phi(e^{i\bar{\vartheta}})\Psi(e^{i\bar{\vartheta}})^{-1}$
where $\bar{\vartheta}\in[0,2\pi)$ is fixed. The function \eq
f_\nu(\overline{e})=\overline{e}-\frac{\nu}{\nu-1}(\overline{e})^{\frac{\nu-1}{\nu}}\eeq
is the infinitesimal contribution at $\bar{\vartheta}$ to the
objective function
$\frac{1}{\nu}\left(\Sc_T^{(1-\nu^{-1})}(\Phi\|\Psi)-\frac{\nu}{\nu-1}\right)$
which is equivalent to $\Sc_T^{(1-\nu^{-1})}(\Phi\|\Psi)$ in
Problem \ref{pb_spectrum_approx}. Note that,
$d\,f_\nu(\overline{e})/ d\,\overline{e}=1-e^{-\frac{1}{\nu}}$
represents the instantaneous rate of change of
$f_\nu(\overline{e})$ at point $\overline{e}$. Moreover,
$d\,f_{\nu_1}(\overline{e})/ d\,\overline{e}\geq
d\,f_{\nu_1}(\overline{e})/ d\,\overline{e}\geq 0$ with $\nu_1\geq
\nu_2$ and $\overline{e}\geq 1$. Accordingly  the larger $\nu$ is,
the more $f_\nu$ penalizes values of $\overline{e}$ greater than
one. Therefore a sufficiently large value of $\nu$ should avoid
solution whose innovation process is greater than $I$ in narrow
ranges of frequencies.

We conclude that $\Phi_\nu$ with $\nu$ small is preferable when
the model for $y$ should be simple in terms of complexity degree
whereas the one with $\nu$ large is preferable when a small
estimation error (also in terms of innovation process) is
required. Regarding the computational complexity, the small $\nu$
is, the better performance is. In fact, to solve equation
(\ref{eq_Newton_step}) it is required to construct a state space
model having state dimension which is proportional to $\nu$, see
\cite[Section VI]{MATRICIAL_ALGTHM_RAMPONI_FERRANTE_PAVON_2009}.

\section{Conclusion}
In this technical note we have presented a new multivariate
divergence family between spectral densities which arises in the
context of optimal prediction. Such a divergence family leads to a
family of solutions to the spectrum approximation problem which
are with bounded {\em McMillan} degree under the mild assumption
that the {\em a priori} spectral density has bounded {\em
McMillan} degree. Finally, a simulation study has been presented
for drawing the application scenarios of this family of spectral
estimators.

\bibliographystyle{plain}
\bibliography{biblio}

\begin{thebibliography}{10}

\bibitem{AMARI_DIFFERENTIAL_GEOMETRIC}
S.~Amari.
\newblock {\em {Differential-Geometrical Methods in Statistics}}.
\newblock Springer-Verlag, Berlin, 1985.

\bibitem{BASU_ROBUST_1998}
A.~Basu, I.~Harris, N.~Hjort, and M.~Jones.
\newblock Robust and efficient estimation by minimising a density power
  divergence.
\newblock {\em Biometrika}, 85(3):549--559, Sep. 1998.

\bibitem{A_NEW_APPROACH_BYRNES_2000}
C.~Byrnes, T.~Georgiou, and A.~Lindquist.
\newblock A new approach to spectral estimation: A tunable high-resolution
  spectral estimator.
\newblock {\em IEEE Trans. Signal Processing}, 48(11):3189--3205, Nov. 2000.

\bibitem{FERRANTE_TIME_AND_SPECTRAL_2012}
A.~{Ferrante}, C.~{Masiero}, and M.~{Pavon}.
\newblock {Time and spectral domain relative entropy: A new approach to
  multivariate spectral estimation}.
\newblock {\em IEEE Trans. Autom. Control}, 57(10):2561--2575, Oct. 2012.

\bibitem{Hellinger_Ferrante_Pavon}
A.~Ferrante, M.~Pavon, and F.~Ramponi.
\newblock Hellinger versus {K}ullback-{L}eibler multivariable spectrum
  approximation.
\newblock {\em IEEE Trans. Autom. Control}, 53(4):954--967, May 2008.

\bibitem{ME_ENHANCEMENT_FERRANTE_2012}
A.~Ferrante, M.~Pavon, and M.~Zorzi.
\newblock A maximum entropy enhancement for a family of high-resolution
  spectral estimators.
\newblock {\em IEEE Trans. Autom. Control}, 57(2):318--329, Feb. 2012.

\bibitem{FERRANTE_STRUCTURED_2012}
A.~Ferrante, M.~Pavon, and M.~Zorzi.
\newblock Structured covariance estimation in high resolution spectral
  analysis.
\newblock In {\em Proc. of Int. Symp. Mathematical Theory of Network and
  Systems, MTNS 2012}. Melbourne, 2012.

\bibitem{GEORGIUO_THE_STRUCUTRE_2002}
T.~Georgiou.
\newblock The structure of state covariances and its relation to the power
  spectrum of the input.
\newblock {\em IEEE Trans. Autom. Control}, 47(7):1056--1066, Jul. 2002.

\bibitem{RELATIVE_ENTROPY_GEORGIOU_2006}
T.~Georgiou.
\newblock {R}elative entropy and the multivariable multidimensional moment
  problem.
\newblock {\em IEEE Trans. Inform. Theory}, 52(3):1052--1066, Mar. 2006.

\bibitem{KL_APPROX_GEORGIUO_LINDQUIST}
T.~Georgiou and A.~Lindquist.
\newblock {K}ullback-{L}eibler approximation of spectral density functions.
\newblock {\em IEEE Trans. Inform. Theory}, 49(11):2910--2917, Nov. 2003.

\bibitem{ITAKURA_SAITO_1968}
F.~Itakura and S.~Saito.
\newblock Analysis synthesis telephony based on the maximum likelihood method.
\newblock In {\em Proceedings of 6th International Congress on Acoustics},
  pages 17--20, Tokyo, Japan, 1968.

\bibitem{GEORGIOU_DISTANCES_PSD_2012}
X.~Jiang, L.~Ning, and T.~Georgiou.
\newblock Distances and riemannian metrics for multivariate spectral densities.
\newblock {\em IEEE Trans. Autom. Control}, 57(7):1723--1735, Jul. 2012.

\bibitem{KAILATH_LIN_EST}
T.~Kailath, A.~Sayed, and B.~Hassibi.
\newblock {\em Linear estimation}.
\newblock Prentice Hall, 2000.

\bibitem{KULLBACK_LEIBLER_1951}
S.~Kullback and R.~Leibler.
\newblock On information and sufficiency.
\newblock {\em The Annals of Mathematical Statistics}, 22(1):79--86, 1951.

\bibitem{LINDQUIST_PREDICTION_ERROR_2007}
A.~Lindquist.
\newblock Prediction-error approximation by convex optimization.
\newblock In A.~Chiuso, A.~Ferrante, and S.~Pinzoni, editors, {\em Modeling,
  Estmation and Control:Festschrift in honor of Giorgio Picci on the occation
  of his sixty-fifth birthday}, pages 265--275. Springer-Verlag, 2007.

\bibitem{LINDQUIST_PICCI}
A.~Lindquist and G.~Picci.
\newblock Linear stochastic systems: A geometric approach to modeling,
  estimation and identification.
\newblock In preparation: preprint available in
  \url{http://www.math.kth.se/~alq/LPbook}.

\bibitem{MATRICIAL_ALGTHM_RAMPONI_FERRANTE_PAVON_2009}
F.~Ramponi, A.~Ferrante, and M.~Pavon.
\newblock A globally convergent matricial algorithm for multivariate spectral
  estimation.
\newblock {\em IEEE Trans. Autom. Control}, 54(10):2376--2388, Oct. 2009.

\bibitem{BETA}
M.~Zorzi.
\newblock {A new family of high-resolution multivariate spectral estimators}.
\newblock {\em IEEE Trans. Autom. Control}, In Press.

\bibitem{ALPHA}
M.~Zorzi.
\newblock Rational approximations of spectral densities based on the {A}lpha
  divergence.
\newblock {\em Math. Control Signals Syst.}, In Press.

\bibitem{ON_THE_ESTIMATION_ZORZI_2012}
M.~Zorzi and A.~Ferrante.
\newblock On the estimation of structured covariance matrices.
\newblock {\em Automatica}, 48(9):2145--2151, Sep. 2012.

\end{thebibliography}
\end{document}